\documentclass[12pt, leqno]{article}
\usepackage{amsmath}
\usepackage{amssymb}

\def\underset#1#2{{\mathrel{\mathop {{}_{} {#2}}\limits_{{#1}_{}}}}}
\def\upplim_#1{\underset{#1}{\overline\lim}\;}
\def\lowlim_#1{\underset{#1}{\underline\lim}\;}
\newcommand{\B}{{\mathbf{B}}}
\newcommand{\C}{{\mathbf{C}}}

\newcommand{\del}{{\partial}}
\newcommand{\delbar}{\bar{\partial}}
\newcommand{\fa}{{\forall}}

\newcommand{\rP}{{\hat{\mathrm{P}}}}

\newcommand{\pnc}{{\mathbf{P}^n(\mathbf{C})}}
\newcommand{\pone}{{\mathbf{P}^1(\mathbf{C})}}

\newcommand{\T}{{\mathbf{T}}}
\newcommand{\Z}{{\mathbf{Z}}}
\newtheorem{cor}[equation]{Corollary}
\newtheorem{dfn}[equation]{\indent {\it Definition}\rm}

\newtheorem{prop}[equation]{\bf Proposition}
\newtheorem{thm}[equation]{\bf Theorem}
\newtheorem{lem}[equation]{\bf Lemma}
\numberwithin{equation}{section}
\newtheorem{rmk}[equation]{\indent {\it Remark}\rm}
\newenvironment{pf}
{{\it Proof.  }\hskip10pt} {\hfill{\it Q.E.D.}\par\vskip+10pt}

\title{Connections and the Second Main Theorem for Holomorphic Curves
\thanks{  Research supported in part by Grant-in-Aid
for Scientific Research (S) 17104001.\hfill\break
\indent 2010 {\it Matheamtical Subject Classification}:
Primary 32H30; Scondary 30D35.
}}
\date{\hfill}
\author{Junjiro Noguchi}
\begin{document}
\parindent12pt
\baselineskip16pt
\maketitle
\begin{abstract}
By means of $C^\infty$-connections we will prove
a general second main theorem and
some special ones for holomorphic curves.
The method gives a geometric proof of H. Cartan's second main theorem
in 1933.
By applying the same method, we will prove some second main theorems
in the case of the product space $(\pone)^2$ of the Riemann sphere.
\end{abstract}
\section{Main Results.}

In this paper we are going to prove a general second main theorem and
some special ones for holomorphic curves.
We begin with the general one. One finds an application of
this method in Y. Tiba \cite{t10}.

{\bf (a)} Let $M$ be a compact complex manifold of dimension $n$ and
let $\T (M)$ denote the
holomorphic tangent bundle over $M$.
Let $\nabla$ be a $C^\infty$ connection in $\T(M)$.
Let $U$ be a domain of the complex plane $\C$.
For a holomorphic curve $f:U \to M$ we have the derivative
(1-jet lift) $f'(z) \in \T(M)_{f(z)}$, and we set inductively
$$
f^{(1)}(z)=f'(z),\quad f^{(k)}(z)=\nabla_{f'(z)} f^{(k-1)}(z),
\quad k=2,3, \ldots .
$$
We define the Wronskian of $f$ with respect to $\nabla$ by
$$
W(\nabla, f)=f^{(1)}(z) \wedge \cdots \wedge f^{(n)}(z) \in K_M^*,
$$
where $K_M^*$ denotes the dual of the canonical bundle $K_M$ over $M$.
Because of its local nature it makes sense to say that
$W(\nabla, f)$ is holomorphic or that $\log|W(\nabla, f)|$ is
subharmonic.

We say that $f$ is $\nabla$-(resp.\ non)degenerate if and only if
$W(\nabla,f)\equiv 0$ (resp.\ $\not\equiv 0$).

Cf.\ \S\S2\&3 for more notation. The first result of this paper is as follows:

\begin{thm}
\label{smt-1}
Let $f:\C \to M$ be a $\nabla$-nondegenerate holomorphic curve and let $D=\sum_i D_i$
be an effective reduced  divisor with only simple normal crossings.
Assume
\begin{enumerate}
\item
$\log|W(\nabla, f)|$ is subharmonic;
\item
every $D_i$ is $\nabla$-totally geodesic.
\end{enumerate}
Then we have
\begin{equation}
\label{smt-2}
T_f(r, L(D))+T_f(r, K_M) \leq \sum_i N_n(r, f^*D_i)+S_f(r).
\end{equation}
\end{thm}

Here $N_n(r, f^*D_i)$ denotes the $n$-truncated counting function of $f^*D_i$,
and $S_f(r)$ a small term in the Nevanlinna
theory such as
$$
S_f(r)=O(\log r + \log T_f(r))||,
$$
with the order function $T_f(r)$ of $f$ with respect an hermitian metric
or an ample line bundle over $M$.

Let $\nabla$ be the Fubini-Study metric connection on the $n$-dimensional
complex projective space $\pnc$. Then $\nabla$-totally geodesic
complex submanifolds of $\pnc$ are complex linear subspaces.
A holomorphic curve $f:\C \to \pnc$ is linearly nondegenerate if and only
if $W(\nabla, f)\not\equiv 0$; moreover, the {\it Wronskian $W(\nabla, f)$
is holomorphic} (see Theorem \ref{w-hol}).

\begin{cor}
\label{pnc}
Let $f:\C \to \pnc$ be a linearly nondegenerate holomorphic curve
and let $1 \leq n \leq 3$.
Then for $q$ hyperplanes $H_i \subset \pnc$, $1 \leq i \leq q$, in
 general position we have
\begin{equation}
\label{cartan}
q T_f(r, O(1))+T_f(r, K_{\pnc}) \leq \sum_i N_n(r, f^*D_i)+S_f(r),
\end{equation}
where $O(1)$ denote the hyperplane bundle over $\pnc$.
\end{cor}

Note that \eqref{cartan} is Cartan's Second Main Theorem
(\cite{c33}), for $K_{\pnc}=O(-n-1)$.
Thus, this gives a {\it geometric proof}
of Cartan's Second Main Theorem.
\medskip

{\bf (b)} It is our second aim to consider a special case where
we deal with a holomorphic curve $f:\C \to \pone^2$ (\S5).
We will consider
$\pone^2$ as an equivariant compactification of the semi-abelian variety
$G={\C^*}^2$. This is a quite special case, but interesting, while
not much study has been done for it in the past.

Let $E=\pone^2 \setminus G$ be the boundary divisor
and let $(x,y)$ be the affine coordinate system of $G$.
Here we denote by $\nabla$ the flat connection with respect to the
invariant vector fields, $x\frac{\del}{\del x}$, $y\frac{\del}{\del y}$
on $G$.

We will prove two Theorems \ref{smt-5} and \ref{smt-7} of
Nevanlinna's second main theorem type for $\nabla$-nondegenerate $f$
with some additional condition
and a divisor $D$ on $\pone^2$ whose irreducible components are all
$\nabla$-totally geodesic.

\begin{rmk}\rm
For a more general case without additional condition, see Y. Tiba
 \cite{t10},
where he applies the same method as above
for a problem to obtain a new second main theorem.
\end{rmk}

At the end we will give some examples and problems.
\smallskip

{\it Acknowledgment.}  The author would like to express his
sincere gratitude to Professor J\"org Winkelmann for
interesting discussions on the present subject.
The author is grateful to the referee for pointing an oversight
in the original proof of Theorem \ref{w-hol},
which was minor but certainly necessary to be fixed.

\section{Totally geodesic divisor and Lemma on logarithmic derivative.}

Let $M$ be a complex $n$-dimensional manifold, and let
$\nabla$ be a $C^\infty$ connection in $\T(M)$; i.e., for $C^\infty$ vector
fields $X, Y$ and a $C^\infty$ function $\alpha$ on $M$ we have
\begin{enumerate}
\item
$\nabla_X Y$ is a $C^\infty$ vector field in $\T(M)$, and is linear
in $X$ and $Y$ over $\C$;
\item
$\nabla_{\alpha X}Y=\alpha \nabla_X Y$;
\item
$\nabla_X(\alpha Y)=X(\alpha)\cdot Y+ \alpha \nabla_X Y$.
\end{enumerate}

Let $N$ be a locally closed complex submanifold of $M$. Take $C^\infty$ sections
$X', Y'$ in $\T(N)$, and extend them to $C^\infty$ sections, $X, Y$
in $\T(M)$ over a neighborhood of $N$ in $M$.
Then the restriction  $(\nabla_X Y)|_N$ is independent of the extensions
$X, Y$, and so denoted by $\nabla_{X'} Y'$, which is a section in $\T(M)|_N$, but
{\it not} in $\T(N)$ generally.
\begin{dfn}
A locally closed complex submanifold $N$ of $M$ is said to be $\nabla$-totally geodesic if
$\nabla_{X'} Y'$ is valued in $\T(N)$ for all $C^\infty$ sections
$X',Y'$ in $\T(N)$.
\end{dfn}

Let $f: \C \to M$ be a holomorphic curve. Then $W(\nabla, f)(z)$ is
valued in the dual $K^*_{M, f(z)}$ of the canonical bundle $K_M$ at $f(z)$.
We take a $C^\infty$ volume form $\Omega$ on $M$.
Then $|W(\nabla, f)(z)|^2 \cdot \Omega(f(z))$ is a non-negative $C^\infty$
function in $z \in \C$.

Let $D=\sum_i D_i$ be a divisor on $M$ with irreducible components $D_i$.
Let $\sigma_i$ be a section of the line bundle
$L(D_i)$ determined by $D_i$ such that the divisor $(\sigma_i)$ coincides
$D_i$, and introduce a hermitian metric $\|\cdot \|$ in every $L(D_i)$.
We set
\begin{equation}
\label{xi}
\xi(z)=\frac{|W(\nabla, f)(z)|^2\cdot \Omega(f(z))}{\prod_i \|\sigma_i(f(z))\|^2}.
\end{equation}
As usual, we set $\log^+ \xi(z)=\log\max \{1, \xi(z)\}$.

The following is a version of Nevanlinna's lemma on logarithmic
derivative (cf.\ \cite{n77}, \cite{no84}):

\begin{lem}
\label{loglem}
Let $M$ be a complex algebraic manifold and let $D=\sum_i D_i$ be a
divisor with irreducible components $D_i$.
Assume that
\begin{enumerate}
\item
$D$ has only simple normal crossings;
\item
every $D_i$ is $\nabla$-totally geodesic.
\end{enumerate}
Then we have
\begin{equation}
\label{lemlog}
\int_{|z|=r}\log^+\xi(z) \frac{d\theta}{2\pi}=S_f(r).
\end{equation}
\end{lem}

\begin{pf}
Let $M=\cup_\alpha U_\alpha$ be a finite affine covering with
rational functions $x_{\alpha}^i, 1 \leq i \leq n=\dim M$ over $M$
such that
\begin{enumerate}
\item
$x_{\alpha}^i, 1 \leq i \leq n$, are holomorphic on $U_\alpha$,
and give rise to coordinates in a neighborhood of every point of
$U_\alpha$;
\item
$U_\alpha \cap D=\{x_\alpha^1 \cdots x_\alpha^{k_\alpha}=0\}$.
\end{enumerate}
Let $V_\alpha \Subset U_\alpha$ be relatively compact open subsets
such that $M=\cup V_\alpha$.
Let $\mathbf{1}_{V_\alpha}$ be the characteristic function of
the set $V_\alpha$.
Set $f_\alpha^j(z)=x_\alpha^j(f(z)), 1 \leq j \leq n$.
Let $\Gamma_{\alpha i j}^k$ be the Christofell symbols of
$\nabla$ with respect to $(x_\alpha^i)$.
Since $D_i$ are $\nabla$-totally geodesic, there are $C^\infty$
functions $A_\alpha$ and $B_\alpha$ on $U_\alpha$ such that
\begin{equation}
\label{geod}
\Gamma_{\alpha i j}^h(x_\alpha^1, \ldots, x_\alpha^h,
\ldots, x_\alpha^n)=
A_{\alpha ij}^h \cdot x_\alpha^h+B_{\alpha ij}^h \cdot \bar{x}_\alpha^h,
\quad 1 \leq h \leq k_\alpha.
\end{equation}
Therefore there is a constant $C_\alpha>0$ such that
\begin{align}
\label{log1}
\left|\frac{\mathbf{1}_{V_\alpha}(f(z))}{f_\alpha^h(z)}\Gamma_{\alpha i j}^h(f(z))\right|
&=\mathbf{1}_{V_\alpha}(f(z))\left|A_{\alpha ij}^h(f(z)) +
B_{\alpha ij}^h(f(z)) \frac{\bar{f}_\alpha^h(z)}{f_\alpha^h(z)}\right|\\
\nonumber
&\leq \mathbf{1}_{V_\alpha}(f(z))\left(|A_{\alpha ij}^h(f(z))| +
|B_{\alpha ij}^h(f(z))|\right)\\
\nonumber
&\leq C_\alpha, \qquad 1 \leq h \leq k_\alpha.
\end{align}
Here we understand ``$\frac{0}{0}=0$'',
when $\mathbf{1}_{V_\alpha}(f(z))=0$
and $f_\alpha^h(z)=0$, provided that $f(z) \in U_\alpha$;
we extend the above function for all $z \in \C$, as zero, when
$f(z) \not\in U_\alpha$.
We set
$$
f^{(l)}(z)= f_\alpha^{(l)k}(z)\left(\frac{\del}{\del x_\alpha^k}
\right)_{f(z)},
$$
where Einstein's convention is used for summation.

There is a $C^\infty$ function $a_\alpha$ on $U_\alpha$ such that
\begin{align}
\label{xi0}
\xi(z)&=\left| \det 
\begin{pmatrix}
f_\alpha^{(1)1} & \cdots & f_\alpha^{(1)k_\alpha}& \cdots &
f_\alpha^{(1)n}\\
f_\alpha^{(2)1} & \cdots & f_\alpha^{(2)k_\alpha}& \cdots &
f_\alpha^{(2)n}\\
f_\alpha^{(3)1} & \cdots & f_\alpha^{(3)k_\alpha}& \cdots &
f_\alpha^{(3)n}\\
\vdots & \vdots & \vdots & \vdots & \vdots \\
f_\alpha^{(n)1} & \cdots & f_\alpha^{(n)k_\alpha}& \cdots &
f_\alpha^{(n)n}
\end{pmatrix}
\right|^2 \frac{a_\alpha(f(z))}{|f_\alpha^1|^2 \cdots
|f_\alpha^{k_\alpha}|^2}\\
\nonumber
&=\left|\left|
\begin{matrix}
\frac{f_\alpha^{(1)1}}{f_\alpha^1} & \cdots & 
\frac{f_\alpha^{(1)k_\alpha}}{f_\alpha^{k_\alpha}}&
f_\alpha^{(1)k_\alpha+1} & \cdots &
f_\alpha^{(1)n}\\
\frac{f_\alpha^{(2)1}}{f_\alpha^1} & \cdots & 
\frac{f_\alpha^{(2)k_\alpha}}{f_\alpha^{k_\alpha}}&
f_\alpha^{(2)k_\alpha+1} & \cdots &
f_\alpha^{(2)n}\\
\frac{f_\alpha^{(3)1}}{f_\alpha^1} & \cdots & 
\frac{f_\alpha^{(3)k_\alpha}}{f_\alpha^{k_\alpha}}&
f_\alpha^{(3)k_\alpha+1} & \cdots &
f_\alpha^{(3)n}\\
\vdots & \vdots & \vdots & \vdots & \vdots & \vdots \\
\frac{f_\alpha^{(n)1}}{f_\alpha^1} & \cdots & 
\frac{f_\alpha^{(n)k_\alpha}}{f_\alpha^{k_\alpha}}&
f_\alpha^{(n)k_\alpha+1} & \cdots &
f_\alpha^{(n)n}\\
\end{matrix}
\right|\right|^2 a_\alpha(f(z)).
\end{align}
Note that $\mathbf{1}_{V_\alpha}(f(z))\cdot f_\alpha^i(z), 1 \leq i \leq n$,
and $\mathbf{1}_{V_\alpha}(f(z)) \cdot a_\alpha(f(z))$ are bounded
functions.
Therefore we have by \eqref{xi0}
\begin{align}
\label{log2}
\log^+\xi(z)&=O\left(\sum_\alpha
\left(\sum_{ 1 \leq k\leq k_\alpha,1\leq l \leq n}
\mathbf{1}_{V_\alpha}(f(z))\cdot \log^+
\left| \frac{f_\alpha^{(l)k}(z)}{f_\alpha^k(z)}\right|
\right. \right. \\
\nonumber
&\quad \left. \left.
+\sum_{1\leq k,l \leq n}\mathbf{1}_{V_\alpha}(f(z))\cdot
\log^+|f_\alpha^{(l)k}(z)|
\right)
\right)+O(1),
\end{align}
where the estimate ``$O(*)$'' is uniform in $z \in \C$
(it is used in the same sense from now on).
We first compute $f_\alpha^{(l)k}$:
For instance, we have
$$
\log^+|f_\alpha^{(1)k}|
=\log^+|f_\alpha^{k}{'}|=\log^+|f_\alpha^{k(1)}|.
$$
For $1 \leq k \leq k_\alpha$ we have
$$
\log^+\left|
\frac{f_\alpha^{(1)k}}{f_\alpha^k}\right|=
\log^+\left|\frac{f_\alpha^{k}{'}}{f_\alpha^k}\right|.
$$

For $l=2$ we have
$$
f_\alpha^{(2)k}={f_\alpha^k}{''}+\Gamma_{\alpha i_1 i_2}^k \circ f\cdot
f_\alpha^{i_1}{'}f_\alpha^{i_2}{'}.
$$
Since $\mathbf{1}_{V_\alpha}\circ f \cdot \Gamma_{\alpha i_1 i_2}
\circ f$
is bounded, we have
\begin{align*}
\mathbf{1}_{V_\alpha}\circ f \cdot |f_\alpha^{(2)k}|
&=\mathbf{1}_{V_\alpha}\circ f \cdot O\left(
|f_\alpha^k{''}|+\left(\sum_{i=1}^n |f_\alpha^{i}{'}|\right)^2
\right),\\
\mathbf{1}_{V_\alpha}\circ f \cdot \log^+ |f_\alpha^{(2)k}|
&=\mathbf{1}_{V_\alpha}\circ f \cdot O\left(
\log^+ |f_\alpha^{k(2)}|+ \sum_{i=1}^n \log^+ |f_\alpha^{i(1)}|
\right)+O(1).
\end{align*}
For $1 \leq k \leq k_\alpha$ we have by \eqref{log1}
\begin{align*}
\mathbf{1}_{V_\alpha}\circ f \cdot \left|\frac{f_\alpha^{(2)k}}{f_\alpha^k}
\right|&=\mathbf{1}_{V_\alpha}\circ f \cdot
O\left(
\left|\frac{f_\alpha^k{''}}{f_\alpha^k}\right|
+\left(\sum_{i=1}^n |f_\alpha^{i}{'}|\right)^2
\right), \\
\mathbf{1}_{V_\alpha}\circ f \cdot \log^+ |f_\alpha^{(2)k}|
&=\mathbf{1}_{V_\alpha}\circ f \cdot O\left(
\log^+ |f_\alpha^{k(2)}|+\sum_{i=1}^n \log^+|f_\alpha^{i(1)}|\right)
+O(1).
\end{align*}

Up to here it is easy to obtain the estimate.
For $l=3$, $f_\alpha^{(3)k}$ starts to involve the partial derivatives of
$\Gamma_{\alpha i_1 i_2}^k$:
\begin{align}
\nonumber
f_\alpha^{(3)k} &= f_\alpha^k{'''}+\Gamma_{\alpha i_1 i_2}^k\circ f\cdot
f_\alpha^{i_1}{''}f_\alpha^{i_2}{'}
+\Gamma_{\alpha i_1 i_2}^k \circ f\cdot
f_\alpha^{i_1}{'}f_\alpha^{i_2}{''}
\\
\nonumber
&\quad +\frac{\del \Gamma_{\alpha i_1 i_2}^k}{\del x_\alpha^{i_3}}
\circ f\cdot
f_\alpha^{i_1}{'}f_\alpha^{i_2}{'}f_\alpha^{i_3}{'}
+\Gamma_{\alpha i_1 i_2}^k \circ f\cdot f_\alpha^{i_1}{'}f_\alpha^{(2)i_2},
\end{align}
and
$$
\Gamma_{\alpha i_1 i_2}^k \circ f\cdot
f_\alpha^{i_1}{'}f_\alpha^{(2)i_2}=
\Gamma_{\alpha i_1 i_2}^k \circ f\cdot
f_\alpha^{i_1}{'}
{f_\alpha^{i_2}}{''}+\Gamma_{\alpha i_1 i_2}^k \circ f\cdot
f_\alpha^{i_1}{'}\cdot
\Gamma_{\alpha i_3 i_4}^{i_2}\circ f\cdot
f_\alpha^{i_3}{'}f_\alpha^{i_4}{'}.
$$
It follows that
\begin{align}
\nonumber
\mathbf{1}_{V_\alpha}\circ f \cdot \log^+ |f_\alpha^{(3)k}|
&=\mathbf{1}_{V_\alpha}\circ f \cdot O\left(
\log^+|f_\alpha^{k(3)}|+\sum_{i=1}^n \log^+ |f_\alpha^{i(2)}|
\right. \\
\nonumber
&\quad \left. +\sum_{i=1}^n\log^+ |f_\alpha^{i(1)}|
\right)+O(1).
\end{align}
For $1 \leq k \leq k_\alpha$ we estimate
$\mathbf{1}_{V_\alpha}\circ f \cdot
\left|\frac{f_\alpha^{(3)k}}{f_\alpha^k}\right|$:
\begin{align}
\label{der3}
&\mathbf{1}_{V_\alpha}\circ f \cdot
\left|\frac{f_\alpha^{(3)k}}{f_\alpha^k}\right|\\
\nonumber
&\leq 
\mathbf{1}_{V_\alpha}\circ f \cdot
\left(\left|\frac{f_\alpha^{k}{'''}}{f_\alpha^k}\right|
+\left|\frac{\Gamma_{\alpha i_1 i_2}^k \circ f}{f_\alpha^k}
f_\alpha^{i_1}{''} f_\alpha^{i_2}{'}\right|
+\left|\frac{\Gamma_{\alpha i_1 i_2}^k  \circ f}{f_\alpha^k}
f_\alpha^{i_1}{'} f_\alpha^{i_2}{''}\right|
\right.
\\
\nonumber
&\quad \left.
+\left|\frac{1}{f_\alpha^k}
\frac{\del \Gamma_{\alpha i_1 i_2}^k}{\del x_\alpha^{i_3}}
\circ f\cdot
f_\alpha^{i_1}{'}f_\alpha^{i_2}{'}f_\alpha^{i_3}{'}\right|
+\left|\frac{\Gamma_{\alpha i_1 i_2}^k  \circ f}{f_\alpha^k}
f_\alpha^{i_1}{'}f_\alpha^{(2)i_2}\right|
\right).
\end{align}
Note that $\mathbf{1}_{V_\alpha}\circ f \cdot
\left|\frac{\Gamma_{\alpha i_1 i_2}^k \circ f}{f_\alpha^k}\right|$ is bounded.
We compute the fourth term in the right hand side of \eqref{der3}:
For $i_3\not=k$ we have
\begin{align*}
&\mathbf{1}_{V_\alpha}\circ f \cdot\left|\frac{1}{f_\alpha^k}
\frac{\del \Gamma_{\alpha i_1 i_2}^k}{\del x_\alpha^{i_3}}
\circ f\cdot
f_\alpha^{i_1}{'}f_\alpha^{i_2}{'}f_\alpha^{i_3}{'}\right|\\
&= \mathbf{1}_{V_\alpha}\circ f \cdot
\left|
\left(\frac{\del A_{\alpha i_1 i_2}^k}{\del x_\alpha^{i_3}}
\circ f
+ \frac{\del B_{\alpha i_1 i_2}^k}{\del x_\alpha^{i_3}}\circ f\cdot
\frac{\bar{f}_\alpha^k}{f_\alpha^k}
\right)
f_\alpha^{i_1}{'}f_\alpha^{i_2}{'}f_\alpha^{i_3}{'}\right|\\
&
\leq
\mathbf{1}_{V_\alpha}\circ f \cdot
\left(
\left|\frac{\del A_{\alpha i_1 i_2}^k}{\del x_\alpha^{i_3}}
\circ f \right|
+ \left|\frac{\del B_{\alpha i_1 i_2}^k}{\del x_\alpha^{i_3}}
\circ f \right|
\right)
|f_\alpha^{i_1}{'}f_\alpha^{i_2}{'}f_\alpha^{i_3}{'}|\\
&= \mathbf{1}_{V_\alpha}\circ f \cdot
O\left(\sum_{i=1}^n |f_\alpha^i{'}|\right)^3.
\end{align*}
For $i_3=k$ we obtain
\begin{align*}
&\mathbf{1}_{V_\alpha}\circ f \cdot\left|\frac{1}{f_\alpha^k}
\frac{\del \Gamma_{\alpha i_1 i_2}^k}{\del x_\alpha^{k}}\circ f\cdot
f_\alpha^{i_1}{'}f_\alpha^{i_2}{'}f_\alpha^{k}{'}\right|\\
&=\mathbf{1}_{V_\alpha}\circ f \cdot\left|
\frac{\del \Gamma_{\alpha i_1 i_2}^k}{\del x_\alpha^{k}}\circ f\cdot
f_\alpha^{i_1}{'}f_\alpha^{i_2}{'}\frac{f_\alpha^{k}{'}}
{f_\alpha^k}\right|\\
&= \mathbf{1}_{V_\alpha}\circ f \cdot
O\left(
\left(\sum_{i=1}^n |f_\alpha^i{'}|\right)^2 
\left|\frac{f_\alpha^k{'}}{f_\alpha^k}\right|
\right).
\end{align*}
Therefore we get
\begin{align}
\nonumber
\mathbf{1}_{V_\alpha}\circ f \cdot \log^+
\left|\frac{f_\alpha^{(3)k}}{f_\alpha^k}\right|
&\leq 
\mathbf{1}_{V_\alpha}\circ f \cdot O
\left( \log^+ \left| \frac{f_\alpha^{k(3)}}{f_\alpha^k}\right|
+ \sum_{i=1}^n \log^+ |f_\alpha^{i(2)}| \right. \\
\nonumber
&\quad \left.
+ \sum_{i=1}^n \log^+ |f_\alpha^{i(1)}|
\right)
+O(1).
\end{align}

In this way we have
\begin{align}
\label{log4}
\mathbf{1}_{V_\alpha}\circ f \cdot \log^+ |f_\alpha^{(l)k}|
&=O\left( \sum_{1 \leq i \leq n, 1 \leq j \leq l}
\mathbf{1}_{V_\alpha}\circ f \cdot \log^+
| f_\alpha^{i(j)}|\right)+O(1),\quad 1 \leq k \leq n,\\
\nonumber
\mathbf{1}_{V_\alpha}\circ f \cdot \log^+ \left|
\frac{f_\alpha^{(l)k}}{f_\alpha^k}\right|
&=O\left( \sum_{1 \leq j \leq l} \mathbf{1}_{V_\alpha}\circ f \cdot \log^+ \left|
\frac{f_\alpha^{k(j)}}{f_\alpha^k}\right| \right.\\
\nonumber
&+ \left. \sum_{1 \leq i \leq n, 1 \leq j \leq l}
\mathbf{1}_{V_\alpha}\circ f \cdot \log^+
| f_\alpha^{i(j)}|\right)+O(1),\quad 1 \leq k \leq k_\alpha.
\end{align}
Notice that for $j \geq 1$
\begin{align}
\label{log5}
\mathbf{1}_{V_\alpha}\circ f \cdot \log^+
| f_\alpha^{i(j)}| &\leq
\mathbf{1}_{V_\alpha}\circ f \cdot \log^+
\left|\frac{f_\alpha^{i(j)}}{f_\alpha^i} \cdot {f_\alpha^i} \right|\\
\nonumber
&= \mathbf{1}_{V_\alpha}\circ f \cdot \log^+
\left|\frac{f_\alpha^{i(j)}}{f_\alpha^i} \right|+O(1).
\end{align}

Combining \eqref{log2}, \eqref{log4} and \eqref{log5}
with Nevanlinna's lemma on logarithmic
derivative (cf., e.g., \cite{no84}), we deduce that
\begin{align*}
\int_{|z|=r} \log^+ \xi(z)\frac{d\theta}{2\pi}
&= O\left(
\sum_{\alpha, 1 \leq k, l \leq n}
\int_{|z|=r}
\mathbf{1}_{V_\alpha}(f(z))\cdot \log^+
\left| \frac{f_\alpha^{k(l)}(z)}{f_\alpha^k(z)}\right|
\frac{d\theta}{2\pi}
\right)+O(1)\\
&= O\left(
\sum_{\alpha, 1 \leq k, l \leq n}
\int_{|z|=r}
\log^+
\left| \frac{f_\alpha^{k(l)}(z)}{f_\alpha^k(z)}\right|
\frac{d\theta}{2\pi}
\right)+O(1)\\
&=S_f(r).
\end{align*}
\end{pf}

\section{Proof of Theorem \ref{smt-1}.}
We first note that the current
$$
dd^c \log |W(\nabla, f)|^2  =\frac{i}{2\pi} \del \delbar  \log
 |W(\nabla, f)|^2
$$
is well defined and is a positive measure on $\C$.
For the sake of notational simplicity we write $c_1(D)$ for the curvature
form of the hermitian line bundle $L(D)$ defining the Chern class.
It follows that
\begin{equation}
\label{3.1}
dd^c\log \xi=f^*c_1(D)+f^*c_1(K_M)-\sum_i f^*D_i +dd^c\log
 |W(\nabla, f)|^2.
\end{equation}

Let $Z=\sum_\nu \lambda_\nu \cdot z_\nu$ be a divisor on $\C$ with
distinct $z_\nu \in \C$. We set the $k$-truncated divisor of $Z$ with
$k \leq \infty$ by
$$
(Z)_k=\sum_\nu \min\{\nu, k\} \cdot z_\nu.
$$

Calculating the multiplicity at $f(z) \in \sum D_i$, we see that
$$
-\sum_i f^*D_i +dd^c\log |W(\nabla, f)|^2\geq -\sum (f^*D_i)_n
$$
as currents.
It follows from this and \eqref{3.1} that
\begin{equation}
\label{3.2}
dd^c\log \xi \geq f^*c_1(D)+f^*c_1(K_M)-\sum_i (f^*D_i)_n.
\end{equation}
We denote by $T_f(r, L(D))$ (resp. $T_f(r, K_M)$)
the order function of $f$ with respect to $c_1(D)$
 (resp. $c_1(K_M)$); e.g.,
$$
T_f(r, L(D))=\int_1^r \frac{dt}{t} \int_{|z|<t} f^* c_1(D).
$$
Using the counting function $N_n(r, f^*D_i)$ truncated to level $n$
defined by
$$
N_n(r, f^*D_i)=\int_1^r \frac{dt}{t}\int_{|z|<t} (f^* D_i)_n,
$$
we have by Jensen's formula
\begin{align}
\label{3.3}
&T_f(r, L(D))+T_f(r, K_M) \\
\nonumber
&\leq \sum_i N_n(r, f^*D_i)
+\frac{1}{2}\int_{|z|=r} \log \xi(z) \frac{d\theta}{2\pi}
-\frac{1}{2}\int_{|z|=1} \log \xi(z) \frac{d\theta}{2\pi}.
\end{align}
By Lemma \ref{loglem} we see that
$$
T_f(r, L(D))+T_f(r, K_M) \leq \sum_i N_n(r, f^*D_i)+S_f(r).
$$
This finishes the proof.

\section{Geometric proof of Cartan's second main theorem}

The purpose of this section is to give a geometric proof
of H. Cartan's second main theorem (\cite{c33}), whose key is
the following.
We let $\nabla$ denote the connection induced by the
Fubini-Study metric form $\omega$ on $\pnc$ in this section. 
\begin{thm}
\label{w-hol}
Let $f:U \to \pnc$ be a holomorphic curve from a domain $U \subset \C$.
\begin{enumerate}
\item
The holomorphic curve $f$ is $\nabla$-degenerate if and only if
$f$ is linearly degenerate; i.e., the image $f(U)$ is contained
in a hyperplane.
\item
The Wronskian $W(\nabla, f)(z)$ is holomorphic in $z \in U$.
\end{enumerate}
\end{thm}

{\it Remark.}
So far by our knowledge, the above (ii) was first proved by Siu \cite{siu87}
in the case of $n=2$. Since there is no reference providing a proof
for general $n$, we give a proof making use of the potential of the
K\"ahler form $\omega$. The statement (i) should be known, but
since we do not know a reference and its proof consists a part
of the proof of (ii), we here give a self-contained proof.

\begin{pf}
(i) We take a point $z_0 \in U$ and set $w_0=f(z_0) \in \pnc$.
We may assume that $z_0=0$.
Let $(w^1, \ldots, w^n)$ be the normalized affine coordinate of $\pnc$
such that $w_0=(0, \ldots, 0)$ and 
$$
\omega=dd^c\log \left(1+ \|w\|^2\right), \quad
\|w\|^2=\sum_{j=1}^n |w^j|^2.
$$

We set
\begin{align*}
\del_z &=\frac{\del}{\del z}, \qquad \delbar_z=\frac{\del}{\del \bar{z}},\\
f(z) &=(f^1(z), \ldots, f^n(z)).
\end{align*}

Setting $\phi(w)=\log \left(1+ \|w\|^2\right)$ and
$\omega= g_{i\bar{j}}\frac{i}{2\pi}dw^i \wedge d\bar{w}^j
$ we have
$$
g_{i\bar{j}}=\del_i \delbar_j \phi(w),
$$
where $\del_i=\del/\del w^i$ and $\delbar_j=\del/\del \bar{w}^j$.
Then the Christoffel symbol $\Gamma_{ij}^k$ of the connection $\nabla$
is given by
\begin{align*}
\Gamma_{ij}^k&= (\del_i g_{j\bar{l}})\cdot g^{\bar{l} k}=
 \del_i g_{j\bar{l}} \cdot g^{\bar{l} k}\: (=\Gamma_{ji}^k),\\
\nabla_{\del_i} \del_j &=\Gamma_{ij} ^k\del_k,
\end{align*}
where  $(g^{\bar{l} k})$ denotes the inverse matrix of $(g_{j\bar{l}})$
with $g_{j\bar{l}}g^{\bar{l}k}=\delta_j^k$ (Kronecker).
For later use we note that
\begin{align}
\label{diffinv}
\del_h g^{\bar{l}k} &=-g^{\bar{l}i}\cdot \del_h g_{i\bar{j}}
\cdot g^{\bar{j}k} ,\\
\nonumber
\delbar_h g^{\bar{l}k} &=-g^{\bar{l}i} \cdot \delbar_h g_{i\bar{j}}
\cdot g^{\bar{j}k} .
\end{align}

Now the power series expansion of the potential $\phi(w)$ about the
origin $0$ is
\begin{equation}
\label{pot}
\phi(w)=\sum_{\mu=1}^\infty \frac{(-1)^{\mu-1}}{\mu}
\left( \sum_i w^i \bar{w}^i \right)^{\mu}.
\end{equation}
From this power series expansion we see that
{\it
\begin{enumerate}
\item[{\rm i)}]
the partial differentiation
of $\phi$ evaluated at $0$,
\begin{equation}
\label{diff1}
\del_{i_1} \cdots \del_{i_p}\delbar_{j_1} \cdots \delbar_{j_q}
\phi(0)
\end{equation}
can be non-zero
only when the partial differentiations $\del_j$ and $\delbar_j$
appear exactly in pairs in the partial differentiation;
\item[{\rm ii)}]
in particular, the values of odd order differentiations of $\phi(w)$ are
all zero at $w=0$.
\item[{\rm iii)}]
if we rewrite \eqref{diff1} as
\begin{equation}
\label{diff2}
(\del_{i_1}\delbar_{i_1})^{\mu_1} \cdots (\del_{i_p}\delbar_{i_p})^{\mu_p}
\phi(0),
\end{equation}
where $i_1, \ldots, i_p$ are distinct, then
the values of \eqref{diff2} are the same for all choices of indices
$i_1, \ldots, i_p$,
because of the symmetry in the variables $w^1, \ldots, w^n$;
\end{enumerate}
}
\noindent
We will use only the above properties i) and ii).
From these we obtain
\begin{align*}
\del_i g_{j\bar{l}}(0)&=\del_i \del_j \delbar_l \phi(0) =0,\\
\Gamma_{ij}^k(0) &=0.
\end{align*}

We write $f^{(j)}(z)=f^{(j)k}(z)\del_k$, locally about $0$.
Note that $f^{(j)k}(z)$ are not holomorphic for $j \geq 2$:
\begin{align}
\label{dercomp1}
f^{(2)k}&=\del_z^2f^k +\del_zf^{i_1}\cdot \del_z f^{i_2}\cdot
\Gamma_{i_1i_2}^k\circ f,\\
\nonumber
f^{(3)k}&=\del_z f^{(2)k} +\del_zf^{i_1}\cdot \del_z f^{(2)i_2}\cdot
\Gamma_{i_1i_2}^k \circ f\\
\nonumber
&= \del_z^3f^k+ 3\del_z^2f^{i_1}\cdot \del_zf^{i_2}\cdot
 \Gamma_{i_1 i_2}^k \circ f\\
\nonumber
&\quad
+ \del_zf^{i_1}\cdot \del_z f^{i_3}\cdot \del_zf^{i_4}\cdot
\Gamma_{i_3 i_4}^{i_2}\circ f \cdot \Gamma_{i_1 i_2}^k \circ f\\
\nonumber
&\quad +\del_z f^{i_1}\cdot \del_z f^{i_2} \cdot \del_z f^\alpha \cdot
\del_\alpha \Gamma_{i_1 i_2}^k \circ f.
\end{align}
Therefore we see that
\begin{align}
\label{dercomp2}
f^{(j)k} &= \del_z f^{(j-1)k}+\del_zf^{i_1}\cdot f^{(j-1)i_2}
\cdot \Gamma_{i_1 i_2}^k \circ f\\
\nonumber
&=\del_z^j f^k + P^{jk}(\del_z f^{h}, \ldots,
\del_z^{j-1}f^h, \del_{\alpha_1} \cdots \del_{\alpha_{\nu}}
\Gamma_{i_1i_2}^h \circ f), \\
\nonumber
&\quad\:  (1 \leq h \leq n,\: 0 \leq \nu \leq j-2),
\end{align}
where $P^{jk}$ is a polynomial such that every term has a factor
of the form $\del_{\alpha_1} \cdots \del_{\alpha_{\nu}}
\Gamma_{i_1i_2}^h \circ f$.
It follows from \eqref{diff1} and \eqref{diffinv} that
$\del_{\alpha_1} \cdots \del_{\alpha_{\nu}}
\Gamma_{i_1i_2}^h(0)=0$.
Therefore,
$$
P^{jk}(\del_z f^{h}(0), \ldots,
\del_z^{j-1}f^h(0), \del_{\alpha_1} \cdots \del_{\alpha_{\nu}}
\Gamma_{i_1i_2}^h(0))=0,
$$
so that
$$
f^{(j)k}(0) =\del_z^j f^k(0).
$$
It follows that $W(\nabla, f)(0) =\det(\del_z^jf^k(0))$.
Therefore $W(\nabla, f)\equiv 0$ if and only if
the standard Wronskian of $f$ $\det(\del_z^j f^k)\equiv 0$,
and hence if and only if $f$ is linearly degenerate.

(ii)
By a unitary transformation of $(w^i)$ we may assume that the matrix
$\left(\del_z^jf^k(0)\right)_{1 \leq j, k, \leq n}$ is of
lower triangle:
\begin{equation}
\label{triang}
\del_z^j f^k (0)=0,\: j<k,\quad  \del_z^k f^k (0)=c_k
\qquad (c_k \in \C).
\end{equation}

We prove
\begin{equation}
\label{dbar}
\delbar_z f^{(j)k}(0)=0,\quad j \leq k \leq n.
\end{equation}
Let $k \geq j$.
Then it follows from \eqref{dercomp2} that
\begin{align*}
\delbar_z f^{(j)k}
&=\tilde{P}^{jk}(\del_z f^{h}, \ldots, \del_z^{j-1}f^h,
\del_{\alpha_1} \cdots \del_{\alpha_{\nu}}\Gamma_{i_1i_2}^h \circ f,\\
&\qquad\qquad
\delbar_z \del_{\alpha_1} \cdots \del_{\alpha_{\nu}}
\Gamma_{i_1i_2}^h \circ f)\\
&= \tilde{P}^{jk}(\del_z f^{h}, \ldots, \del_z^{j-1}f^h,
\del_{\alpha_1} \cdots \del_{\alpha_{\nu}}\Gamma_{i_1i_2}^h \circ f,\\
&\qquad\qquad
\overline{\del_z f^{\beta}} \cdot
\delbar_{\beta}\del_{\alpha_1} \cdots \del_{\alpha_{\nu}}
\Gamma_{i_1i_2}^h \circ f),
\end{align*}
where $\tilde{P}^{jk}$ is a polynomial naturally derived from
$P^{jk}$ by differentiations.
Therefore we see that
\begin{align}
\label{dbar0}
\delbar_z f^{(j)k}(0) 
= \tilde{P}^{jk}( &\del_z f^{h}(0), \ldots,
\del_z^{j-1}f^h(0),\del_{\alpha_1} \cdots \del_{\alpha_{\nu}}
\Gamma_{i_1i_2}^h(0) , \\
\nonumber
& \overline{\del_z f^{\beta}}(0) \cdot
\delbar_{\beta}\del_{\alpha_1} \cdots \del_{\alpha_{\nu}}
\Gamma_{i_1i_2}^h(0)).
\end{align}
One infers from this, \eqref{diffinv}, \eqref{diff1},
and \eqref{triang} that
the remaining terms in \eqref{dbar0}
are only those involving
$$
\del_z^lf^{i_2}(0)\cdot \overline{\del_z f^{1}}(0) \cdot
\delbar_{1}\Gamma_{1 i_2}^k(0), \quad i_2 \leq l \leq j-1
\; (<k). \hbox{ (Cf.\ \eqref{dercomp1}.)}
$$
Since
$\delbar_{1}\Gamma_{1 i_2}^k(0)=\delbar_1\del_1 \del_{i_2}\delbar_k\phi(0)=0$
for $i_2\not=k$, we have proved \eqref{dbar}.

We finally see that
$$
\delbar_zW(\nabla, f)(0)
= \sum_{j=1}^n
\begin{vmatrix}
c_1 & 0 & \empty & \cdots & \empty & \empty & 0 \\
* & \ddots & \ddots & \empty & \empty & \empty & \empty \\
\empty & \empty & c_{j-1} & 0 & \empty & \empty \\
\vdots & \empty & \empty & 0 & 0 & \empty & \vdots \\
\empty & \empty & \empty & \empty & c_{j+1} & \ddots & \empty \\
\empty & \empty & \empty & \empty & \empty & \ddots & 0 \\
* & \empty & \empty & \cdots & \empty & * & c_n
\end{vmatrix}
=0
$$
\end{pf}

{\it Remark.} Let $\B=\{\|x\|<1\}$ be the unit ball of $\C^n$ with the Bergman
metric $(h_{i\bar{j}})$ on $\B$. Then we have
\begin{align}
\label{potb}
\psi &= \log(1-\|x\|^2)=\sum_{\nu=1}^\infty \frac{1}{\nu}\|x\|^{2\nu},\\
\nonumber
h_{i\bar{j}} &=\del_i \delbar_j \psi .
\end{align}
Let $\nabla_{\B}$ be the connection induced from the Bergman metric
on $\B$. Because of the type of the power expansions \eqref{pot}
and \eqref{potb}, we have
\begin{cor}
The Wronskian $W(\nabla_{\B}, f)$ is holomorphic
for a holomorphic curve $f:U \to \B$.
\end{cor}

\section{Holomorphic curves into $\pone^2$.}

In this section we set
$$
G=\C^{*2},
$$
which is a two-dimensional semi-abelian variety.
We consider $\pone^2$ as an equivariant compactification of $G$.
We fix an affine coordinate system $(x,y)\in G \subset \pone^2$.
Then there are invariant vector fields on $G$,
$$
X=x\frac{\del}{\del x}, \quad Y=y\frac{\del}{\del y},
$$
which form a frame of the holomorphic tangent bundle $\T(G)$.
In this section, we denote by $\nabla$ the flat connection with respect
to the frame $\{X, Y\}$ of $\T(G)$; i.e.,
$$
\nabla_XY = \nabla_Y X=0.
$$
Then $\nabla$ is a meromorphic connection with logarithmic poles
along the boundary divisor $\del G=E$,
which has only simple normal crossings.

We set locally
$$
u=\log x, \qquad v=\log y.
$$
A locally closed complex submanifold $N$ of $G$ is $\nabla$-totally
geodesic if and only if $N$ is an open subset of an affine
linear subspace $\{(u,v)\in \C^2; \lambda u+ \mu v=c\}$
with constants $\lambda, \mu$, and $c$; in particular, $N$
is an open subset of a translate of an analytic 1-parameter
subgroup of $G$.

Let $D \subset G$ be an algebraic reduced divisor, and denote
by the same $D$ the closure in $\pone^2$.
We are going to deal with the Nevanlinna theory for an algebraically
nondegenerate holomorphic curve $f: \C \to \pone$ and for $D+E$;
in particular, we are interested in the problem of the possible
second main theorem.

If $f(\C)\cap E=\emptyset$, then we have $f:\C \to G$.
In this case we know the following theorem by \cite{nwy02} and
\cite{nwy08}.

\begin{thm}
\label{smt-sa}
Assume that $f:\C \to G$ is algebraically nondegenerate,
and let $D$ be an algebraic reduced divisor on $G$.
Then there an equivariant compactification $\hat{G}$ of $G$
such that
$$
T_{\hat{f}}(r, L(\hat{D})) \leq N_1(r, f^*D)+
\epsilon T_{\hat{f}}(r, L(\hat{D})) ||_\epsilon,\quad
\fa \epsilon >0,
$$
where $\hat{f}=f:\C \to \hat{G}$ and
$\hat{D}$ is the closure of $D\cap G$ in $\hat{G}$.
\end{thm}

As the second main theorem for $f: \C \to G$,
Theorem \ref{smt-sa} is the best possible result.
Therefore in the sequel we will be mainly interested in the case where
$
f(\C) \cap E \not= \emptyset.
$
We set
$$
f(z)=(F(z), G(z))
$$
with respect to the coordinate system $(x,y)$,
where $F(z)$ and $G(z)$ are meromorphic functions in $\C$.
It follows that
\begin{equation}
\label{nablawr}
W(\nabla, f)=
\begin{vmatrix}
\frac{F'}{F} & \frac{G'}{G} \\
\left(\frac{F'}{F}\right)' & \left(\frac{G'}{G}\right)'
\end{vmatrix}
\cdot F \cdot G \cdot
\left( \frac{\del}{\del x} \wedge \frac{\del}{\del y}
\right)_{f(z)}.
\end{equation}

\begin{prop}
\label{deg}
If $f$ is $\nabla$-degenerate and $f(\C)\cap E \not= \emptyset$,
then $f(\C)$ is contained in the closure of a translate of
a $1$-dimensional algebraic subgroup of $G$.
\end{prop}

\begin{pf}
Suppose that $W(\nabla, f)\equiv 0$.
Then there is a non-trivial linear relation with
$\lambda , \mu \in \C$:
\begin{equation}
\label{rel-1}
\lambda \frac{F'(z)}{F(z)}+ \mu \frac{G'(z)}{G(z)}=0, \quad z \in \C.
\end{equation}
If one of $\lambda$ and $\mu$ is zero, the conclusion is immediate.
Thus we assume that $\lambda \mu \not=0$.
Since $f(\C)\cap E \not=\emptyset$, there is a point $a \in \C$
with $f(a) \in E$. Then $F(a)=0$, or $\infty$, or $G(a)=0$, or $\infty$.
Assume that $F(a)=0$; the other cases are dealt similarly.
Then there are an integer $m$ and a non-vanishing holomorphic function
$\tilde{F}$ in a neighborhood of $a$ such that $F(z)=(z-a)^m \tilde{F}(z)$,
locally.
It follows from \eqref{rel-1} that
\begin{equation}
\label{rel-2}
\frac{\lambda m}{z-a}+\mu \frac{G'(z)}{G(z)}
\end{equation}
is holomorphic about $a$. Therefore $\frac{G'(z)}{G(z)}$ must have a
 pole at $a$. Hence there are a non-zero integer $n$ and a non-vanishing
holomorphic function $\tilde{G}$ in a neighborhood of $a$ with
$G(z)=(z-a)^n \tilde{G}(z)$.
We infer from \eqref{rel-2} that
$$
\frac{\lambda m}{z-a}+\frac{\mu n}{z-a}
$$
is holomorphic about $a$, so that
$$
\lambda m + \mu n=0.
$$
Combining this with \eqref{rel-1}, we have
$$
F(z)^m G(z)^n=c, \quad z \in \C,
$$
where $c \in \C^*$ is a constant.
Thus $f(\C)$ is contained in the closure of a translate of the
algebraic subgroup $\{x^m y^n=1\}$ of $G$.
\end{pf}

Let $D=\sum D_i$ be a divisor on $G$ with only simple normal crossings,
where $D_i$ are the irreducible components.
We assume that {\it every $D_i$ is $\nabla$-totally geodesic}. 
We deal with the value distribution of $f$ for $D+E$ in two ways.

{\bf (1)} Here we use an equivariant blow-up of the compactification
$\pone^2$ of $G$.
By \cite{nwy08} (or directly in this case) we have the following:
\begin{lem}
\label{blowup}
Let the notation be as above. Then there is an equivariant
blow-up $\pi:\hat{G} \to \pone^2$ such that $\hat{D}+\hat{E}$
has only simple normal crossings, where $\hat{D}$ is the closure
of $D$ in $\hat{G}$ and $\hat{E}=\hat{G}\setminus G$.
Moreover, if the stabilizer $\{a \in G; a+D=D\}$ of $D$
is finite, then $\hat{D}$ is ample on $\hat{G}$.
\end{lem}

Let $f:\C \to \pone^2$ be a holomorphic curve such that
$f(\C) \not\subset E$.
Then there is a lifting $\hat{f}:\C \to \hat{G}$ with
$\hat{G}$ in Lemma \ref{blowup} such that $f=\pi \circ \hat{f}$.
Let $\hat{E}=\sum_{j} \hat{E}_j$ be the irreducible decomposition,
and denote by $\{\rP_k\}$ all the crossing points of the $\hat{E}_i$'s.

\begin{thm}
\label{smt-5}
Let $\hat{f}:\C \to \hat{G}$, $\hat{D}$, $\hat{E}$ and $\{\rP_k\}$
be as above.
Assume that $f$ is $\nabla$-nondegenerate.
Then we have
\begin{align}
\nonumber
T_{\hat{f}}(r, L(\hat{D})) &\leq
\sum_{i} N_2(r, \hat{f}^* \hat{D_i})+2 \sum_{j}N_1(r, \hat{f}^*\hat{E}_j)\\
\nonumber
&\quad -\sum_k N_1(r, \hat{f}^* \rP_k)+S_f(r).
\end{align}
\end{thm}

Since $K_{\hat{G}}=-\hat{E}$ and
$$
T_{\hat{f}}(r, K_{\hat{G}}) \leq - \sum_{j}N_1(r, \hat{f}^*\hat{E}_j),
$$
we have the following, formulated closer to the fundamental conjecture
for holomorphic curves (\cite{n04}, \S2 and \cite{n08}).
\begin{cor}
\label{smt-6}
Let the notation be as in Theorem \ref{smt-5}. Then
\begin{equation}
\label{smt-6.1}
T_{\hat{f}}(r, L(\hat{D}))+ 2T_{\hat{f}}(r, K_{\hat{G}})
\leq
\sum_{i} N_2(r, \hat{f}^* \hat{D_i})
+S_f(r).
\end{equation}
\end{cor}

{\it Remark.} (i) The coefficient ``2'' in \eqref{smt-6.1}
should be ``1'' by the fundamental conjecture.

(ii) By Proposition \ref{deg} for $f$ to be $\nabla$-nondegenerate
it suffices to assume
that $f$ is algebraically nondegenerate and $f(\C)\cap E \not= \emptyset$.

For the proof of Theorem \ref{smt-5} we take a holomorphic section
$\hat{\sigma} \in H^0(\hat{G}, L(\hat{D}))$ and
$\hat{\tau}\in H^0(\hat{G}, \hat{E})$ such that the divisors 
$(\hat{\sigma})=\hat{D}$ and $(\hat{\tau})=\hat{E}$.
We introduce a hermitian metric $\|\hat{\sigma}\|$
(resp.\ $\|\hat{\tau}\|$) in $L(\hat{D})$ (resp.\ $L(\hat{E})$).
We take a $C^\infty$ volume form $\Omega$ on $\hat{G}$.
Then we set
\begin{equation}
\label{xi-0}
\hat{\xi}(z)=\frac{|W(\nabla, \hat{f})(z)|^2\cdot \Omega(f(z))}
{\|\hat{\sigma}(f(z))\|^2 \cdot \|\hat{\tau}(f(z))\|^2}.
\end{equation}

\begin{lem}
Let $\hat{\xi}$ be as above \eqref{xi-0}. Then
$$
\int_{|z|=r} \log^+ \hat{\xi}(z) \frac{d\theta}{2\pi}=S_{\hat{f}}(r).
$$
\end{lem}

\begin{pf}
Note that the singularities of $\hat{\xi}(z)$ are those coming from
the intersections of $\hat{f}$ and $\hat{D}+\hat{E}$.
Locally on $\hat{G}$ with local coordinate $\hat{x}, \hat{y}$
such that $\hat{E}$ is written by $\hat{x}=0$, by
$\hat{y}=0$ or by $\hat{x}\hat{y}=0$,
those singularities are given by
$$
\left|
\begin{matrix}
\frac{\frac{d}{dz}\hat{x}(f(z))}{\hat{x}(f(z))}&
\frac{\frac{d}{dz}\hat{y}(f(z))}{\hat{y}(f(z))}\\
\frac{d}{dz}\left(\frac{\frac{d}{dz}\hat{x}(f(z))}{\hat{x}(f(z))}\right)&
\frac{d}{dz}\left(\frac{\frac{d}{dz}\hat{y}(f(z))}{\hat{y}(f(z))}\right)
\end{matrix}
\right|, \quad 
\left|
\begin{matrix}
\frac{\frac{d}{dz}\tilde{\hat{\sigma}}(f(z))}{\tilde{\hat{\sigma}}(f(z))}&
\frac{\frac{d}{dz}\hat{y}(f(z))}{\hat{y}(f(z))}\\
\frac{\frac{d^2}{dz^2}\tilde{\hat{\sigma}}(f(z))}
{\tilde{\hat{\sigma}}(f(z))}&
\frac{d}{dz}\left(\frac{\frac{d}{dz}\hat{y}(f(z))}{\hat{y}(f(z))}\right)
\end{matrix}
\right|,
$$
or by
$$
\left|
\begin{matrix}
\frac{\frac{d}{dz}\tilde{\hat{\sigma}}(f(z))}{\tilde{\hat{\sigma}}(f(z))}&
\frac{\frac{d}{dz}\hat{x}(f(z))}{\hat{y}(f(z))}\\
\frac{\frac{d^2}{dz^2}\tilde{\hat{\sigma}}(f(z))}
{\tilde{\hat{\sigma}}(f(z))}&
\frac{d}{dz}\left(\frac{\frac{d}{dz}\hat{x}(f(z))}{\hat{x}(f(z))}\right)
\end{matrix}
\right|,
$$
where $\tilde{\hat{\sigma}}$ is the local expression of
${\hat{\sigma}}$.
Therefore, from Nevanlinna's Lemma on logarithmic derivative we  deduce
the required estimate
(cf.\ the proof of Lemma \ref{loglem}).
\end{pf}

{\it Proof of Theorem \ref{smt-5}}.
By a careful computation of pole orders in $\hat{\xi}$ we have
the following current inequality on $\C$:
$$
dd^c\log \hat{\xi} \geq \hat{f}^*c_1(\hat{D}) -\sum_i (\hat{f}^*D_i)_2
-2\sum_{j} (\hat{f}^*E_j)_1 + \sum_k (\hat{f}^*\rP_k)_1.
$$
By making use of Jensen's formula and Lemma \ref{xi-0} we complete
the proof of the present theorem. \hfill {\it Q.E.D.}
\smallskip

{\bf (2)}  The advantage of Theorem \ref{smt-5} is that it is applicable
for an arbitrary algebraically non-degenerate $f:\C \to \pone^2$
with $f(\C) \cap E\not=\emptyset$ (cf.\ Proposition \ref{deg}).
On the other hand it is not so easy to compute the order function
$T_{\hat{f}}(r, L(\hat{D}))$.
The blow-up $\hat{G}$ was used to get a kind of the general position
condition with respect $\hat{D}$ and $\hat{f}$.
In the present subsection we are going to deal with the problem
imposing such a condition for $f:\C \to \pone^2$ relative to $D$,
without using a blow-up.

Let $E=\pone^2 \setminus G$ be the boundary divisor of $G$, which
has four components $E_j$ $(1 \leq j \leq 4)$ with
only simple normal crossings at four points,
$$
\mathrm{P}_1=(0,0),\:
\mathrm{P}_2=(0,\infty),\: \mathrm{P}_3=(\infty,0),\: \mathrm{P}_4=(\infty,\infty).
$$

Let $O(m,n)$ denote the line bundle of degree $m$ in the first factor
of $\pone^2$ and of degree $n$ in the second factor of $\pone^2$.
Then every $E_j$ has bidegree $(1,0)$ or $(0,1)$.
Let $\sigma \in H^0(\pone^2, O(m,n))$ such that the divisor $D=(\sigma)$
is reduced and has no common component with $E$.
We introduce the natural metric $\|\sigma\|$.
Let $D=\sum_i D_i$ be the irreducible decomposition of $D$.

\begin{thm}
\label{smt-7}
Let $D=\sum_i D_i$ be as above, and let $f=(F,G):\C \to \pone^2$
be a $\nabla$-nondegenerate holomorphic curve.
We assume the following conditions for $D$ and $f$:
\begin{enumerate}
\item
$D \cap G$ has only simple normal crossings.
\item
Every $D_i$ (strictly speaking, $D_i \cap G$) is $\nabla$-totally
geodesic.
\item
There is a neighborhood $V$ of $\{\mathrm{P}_k\}_{k=1}^4$
such that $f(\C) \cap V =\emptyset$.
\end{enumerate}
Then
\begin{equation}
\label{smt-8}
T_{f}(r, O(m,n)) \leq
\sum_{i} N_2(r, f^* D_i)+2 \sum_{j}N_1(r, f^* E_j)
+S_f(r).
\end{equation}
\end{thm}

\begin{cor}
\label{smt-9}
Under  the same conditions as in Theorem \ref{smt-7} we have,
in particular,
$$
T_{f}(r, O(m-4,n-4)) \leq
\sum_{i} N_2(r, f^* D_i)+S_f(r).
$$
\end{cor}

Let $\Omega$ be the volume form associated with the product of
the Fubini-Study K\"ahler form on $\pone$:
$$
\Omega=\frac{\left(\frac{i}{2\pi}\right)^2dx \wedge d\bar{x}\wedge dy
\wedge d\bar{y}}{\left(1+|x|^2\right)^2\left(1+|y|^2\right)^2}.
$$
Let $\tau \in H^0(\pone^2, O(2,2))$ such that $(\tau)=E$, and
set
\begin{equation}
\label{xi-1}
\xi(z)=\frac{|W(\nabla, f)(z)|^2\cdot \Omega(f(z))}
{\|\sigma(f(z))\|^2 \cdot \|\tau(f(z))\|^2}.
\end{equation}

As in (1), Theorem \ref{smt-7} follows from the following lemma.
\begin{lem}
\label{loglem-2}
Let $\xi(z)$ be as above in \eqref{xi-1}. Then we have
$$
\int_{|z|=r}\log^+\xi(z) \frac{d\theta}{2\pi}=S_f(r).
$$
\end{lem}
\begin{pf}
Notice that $D \cap E=\{\mathrm{P}_j\}_{j=1}^4$.
Let $F$ be the set of intersection points of the irreducible
components of $D\cap G$.
Setting
$$
V_1 = \{\delta < |x| < \delta^{-1}\}\times
\{\delta < |x| < \delta^{-1}\}
$$
with $0 < \delta <1$, 
we take and fix a small $\delta$ so that
\begin{equation}
\label{delta}
F \subset V_1, \qquad
\overline{f(\C)} \cap D \subset V_1.
\end{equation}
In a neighborhood $U (\subset V_1)$ of every point of $F$, $D \cap U$
is defined by
$$
x^{m_i}y^{n_i}=c_i, \quad i=1,2,\quad
|A| \not= 0,
$$
where
$m_i, n_i \in \Z$, $c_i \in \C^*$ and
$A =\begin{pmatrix}
m_1 & m_2 \\
n_1 & n_2
\end{pmatrix}
$.

We first estimate $\xi(z)$, provided $f(z) \in U$.
It follows from \eqref{xi-1} that
with a positive $C^\infty$ function $b$ on $U$
\begin{align}
\label{xi-2}
\xi(z) &= \left|\left|
\begin{matrix}
\frac{F'}{F} & \frac{G'}{G} \\
\left(\frac{F'}{F}\right)' & \left(\frac{G'}{G}\right)'
\end{matrix}
\right|\right|^2
\cdot \frac{b(f(z))}
{|F^{m_1}G^{n_1} - c_1|^2\cdot |F^{m_2}G^{n_2} - c_2|^2}\\
&=||A||^{-2} \left|\left|
\begin{matrix}
m_1 \frac{F'}{F}+n_1 \frac{G'}{G} & m_2 \frac{F'}{F}+n_2 \frac{G'}{G} \\
\nonumber
\left(m_1 \frac{F'}{F}+n_1 \frac{G'}{G}\right)' & 
\left(m_2 \frac{F'}{F}+n_2 \frac{G'}{G}\right)'
\end{matrix}\right|\right|^2\\
\nonumber
&\quad \times \frac{b(f(z))}
{|F^{m_1}G^{n_1} - c_1|^2\cdot |F^{m_2}G^{n_2} - c_2|^2}\\
\nonumber
&=||A||^{-2} \left|\left|
\begin{matrix}
\frac{(F^{m_1}G^{n_1})'}{F^{m_1}G^{n_1}} &
 \frac{(F^{m_2}G^{n_2})'}{F^{m_2}G^{n_2}} \\
\frac{(F^{m_1}G^{n_1})^{''}}{F^{m_1}G^{n_1}}-
\left(\frac{(F^{m_1}G^{n_1})'}{F^{m_1}G^{n_1}}\right)^2
&
\frac{(F^{m_2}G^{n_2})^{''}}{F^{m_2}G^{n_2}}-
\left(\frac{(F^{m_2}G^{n_2})'}{F^{m_2}G^{n_2}}\right)^2 &
\end{matrix}
\right|\right|^2\\
\nonumber
&\quad \times 
\frac{b(f(z))}
{|F^{m_1}G^{n_1} - c_1|^2\cdot |F^{m_2}G^{n_2} - c_2|^2}
\\
\nonumber
&=|| A||^{-2}
\left|\left|
\begin{matrix}
\frac{(F^{m_1}G^{n_1}-c_1)'}{F^{m_1}G^{n_1}-c_1}  \\
\frac{(F^{m_1}G^{n_1}-c_1)^{''}}{F^{m_1}G^{n_1}-c_1}-
\frac{(F^{m_1}G^{n_1})'}{F^{m_1}G^{n_1}}
\cdot \frac{(F^{m_1}G^{n_1}-c_1)'}{F^{m_1}G^{n_1}-c_1}
\end{matrix}
\right.\right. \\
\nonumber
&\left.\left.
\qquad\qquad
\begin{matrix}
\frac{(F^{m_2}G^{n_2}-c_2)'}{F^{m_2}G^{n_2}-c_2} \\
\frac{(F^{m_2}G^{n_2}-c_2)^{''}}{F^{m_2}G^{n_2}-c_2}-
\frac{(F^{m_2}G^{n_2})'}{F^{m_2}G^{n_2}}
\cdot \frac{(F^{m_2}G^{n_2}-c_2)'}{F^{m_2}G^{n_2}-c_2}
\end{matrix}
\right|\right|^2 \\
\nonumber
&\quad \times |F|^{-m_1-m_2} |G|^{-n_1-n_2} b(f(z)).
\end{align}
Since $|F(z)|^{\pm 1}$ and $|G(z)|^{\pm 1}$ are uniformly bounded from above
 by $\delta^{-1}$, provided $f(z)\in V_1$,
we have by \eqref{xi-2}
\begin{align}
\label{xi-3}
\xi(z)  &\leq  P_1 \left(
\left|\frac{(F^{m_1}G^{n_1}-c_1)'}{F^{m_1}G^{n_1}-c_1} \right| ,
\left|\frac{(F^{m_2}G^{n_2}-c_2)'}{F^{m_2}G^{n_2}-c_2} \right| ,
 \right. \\
\nonumber
& \left.
\qquad\quad
\left|\frac{(F^{m_1}G^{n_1}-c_1)^{''}}{F^{m_1}G^{n_1}-c_1}\right| ,
\left|\frac{(F^{m_2}G^{n_2}-c_2)^{''}}{F^{m_2}G^{n_2}-c_2}\right|
\right),
\end{align}
provided $f(z) \in U$,
where $P_1(\cdots)$ is a polynomial with positive coefficients.
Hence, we may assume that \eqref{xi-3} holds, provided
$f(z) \in V$.

In a neighborhood $U'$ of a point of $D \cap (V_1 \setminus V)$,
there is only one irreducible component of $D \cap U'$,
to say, given by
$$
x^{m_1}y^{n_1}=c_1, \quad m_1\not=0.
$$
Then we have
$$
\xi(z)=\frac{1}{m_1^2}
\left|\left|
\begin{matrix}
\frac{(F^{m_1}G^{n_1}-c_1)'}{F^{m_1}G^{n_1}-c_1} &
\frac{G'}{G} \\
\frac{(F^{m_1}G^{n_1}-c_1)^{''}}{F^{m_1}G^{n_1}-c_1}-
\frac{(F^{m_1}G^{n_1})'}{F^{m_1}G^{n_1}}
\cdot \frac{(F^{m_1}G^{n_1}-c_1)'}{F^{m_1}G^{n_1}-c_1} &
\left(\frac{G'}{G}\right)'
\end{matrix}
\right|\right|^2
\cdot |F|^{-2m_1} b(f(z)),
$$
provided $f(z) \in U'$.
Therefore we see that \eqref{xi-3} holds for $f(z) \in V_1$.

Set
\begin{align*}
V_2 &=\{0 \leq |x|\leq \delta\}\times\{\delta \leq |y| \leq \delta^{-1}\}
 \cup
\{\delta \leq |x|\leq \delta^{-1}\}\times\{0 \leq |y| \leq \delta\}\\
&\quad \cup
\{\delta^{-1} \leq |x|\leq \infty \}\times\{\delta \leq |y| \leq
\delta^{-1}\}
 \cup
\{\delta \leq |x|\leq \delta^{-1}\}\times\{\delta^{-1}
 \leq |y| \leq \infty\}.
\end{align*}
It follows from the condition that
$$
f(\C) \subset V_1 \cup V_2.
$$
Note that there exists  a positive constant $c_3$ such that
$$
\|\sigma\|\geq c_3 \hbox{ on } V_2.
$$
Suppose that $f(z) \in V_2$.
Then we obtain
\begin{align*}
\label{xi-4}
\xi(z) &=
\left|\left|
\begin{matrix}
\frac{F'}{F} & \frac{G'}{G}\\
\left(\frac{F'}{F}\right)' &
\left(\frac{G'}{G}\right)'
\end{matrix}
 \right|\right|^2
\cdot \frac{1}{\|\sigma(f(z))\|^2}\\
 & \leq \frac{1}{c_3^2}
\left|\left|
\begin{matrix}
\frac{F'}{F} & \frac{G'}{G}\\
\left(\frac{F'}{F}\right)' &
\left(\frac{G'}{G}\right)'
\end{matrix}
 \right|\right|^2.
\end{align*}
We see by this that there is a polynomial
 $P_2(\cdots)$ with positive coefficients satisfying
\begin{equation}
\label{xi-5}
\xi(z) \leq P_2 \left(
\left|
\frac{F'}{F}
\right|, 
\left|\frac{G'}{G}\right|, 
 \left|
\frac{F^{''}}{F}
\right|,
\left|\frac{G^{''}}{G}\right|
\right).
\end{equation}

Set $P= P_1+P_2$.
We see by \eqref{xi-3} and \eqref{xi-5} that
\begin{align}
\label{xi-6}
\xi(z) &\leq P \left(
\left|
\frac{(F^{m_1}G^{n_1}-c_1)'}{F^{m_1}G^{n_1}-c_1}
\right|,
\left|
\frac{(F^{m_2}G^{n_2}-c_2)'}{F^{m_2}G^{n_2}-c_2}
\right|, \left|\frac{F'}{F}\right|, 
\left|\frac{G'}{G}\right|, 
\right. \\
\nonumber
& \qquad\quad \left.
\left|
\frac{(F^{m_1}G^{n_1}-c_1)^{''}}{F^{m_1}G^{n_1}-c_1}
\right|,
 \left|
\frac{(F^{m_2}G^{n_2}-c_2)^{''}}{F^{m_2}G^{n_2}-c_2}
\right|,
\left|\frac{F^{''}}{F}\right|, 
\left|\frac{G^{''}}{G}\right|
\right)
\end{align}
for all $z \in \C$.
Applying Nevanlinna's lemma on logarithmic derivatives, we infer that
\begin{align*}
&\int_{|z|=r} \log^+\xi(z)\frac{d\theta}{2\pi}\\
& =O\left(
m\left(r, \frac{(F^{m_1}G^{n_1}-c_1)'}{F^{m_1}G^{n_1}-c_1}\right)
+m\left(r,
\frac{(F^{m_2}G^{n_2}-c_2)'}{F^{m_2}G^{n_2}-c_2}\right)
\right. \\
& \quad  +m \left(r, \frac{F'}{F}\right)+
m\left(r, \frac{G'}{G}\right)
 \\
&\quad + m\left( r,
\frac{(F^{m_1}G^{n_1}-c_1)^{''}}{F^{m_1}G^{n_1}-c_1}
\right)+
m \left(r,
\frac{(F^{m_2}G^{n_2}-c_2)^{''}}{F^{m_2}G^{n_2}-c_2}
\right)\\
& \quad \left. 
+m\left(r, \frac{F^{''}}{F}\right)+
m\left(r, \frac{G^{''}}{G}\right)
\right)\\
& =S_f(r).
\end{align*}
\end{pf}
\smallskip

{\it Example 1}.  Let $f(z)=(F(z), G(z))$ be defined by
$$
F(z)=e^z, \qquad G(z)=\frac{e^{z}+1}{e^z-1}.
$$
It is easy to check that $f$ is $\nabla$-nondegenerate,
and to see that the image $f(\C)$ is contained by a curve
$C \subset \pone^2$ defined by
$$
(x-1)(y-1)=2.
$$
Since $C \cap \{\mathrm{P}_j\}_{j=1}^4=\emptyset$, $f$ satisfies the conditions
of Theorem \ref{smt-7}, although $f$ is algebraically degenerate.

{\it Example 2}.  Here we give an example of an algebraically
nondegenerate $f: \C \to \pone^2$ for Theorem \ref{smt-7}.
By Fatou's example (\cite{f22}) there is an injective holomorphic
map $\Phi:\C^2 \to \C^2$ with non-empty exterior open subset
$U$ ($\subset \C^2$) of the image $\Phi(\C^2)$.
By making use of Picard's theorem (or Casorati-Weierstrass' theorem)
we take four points, $(a_i, b_j)\in U$ ($i,j=1,2$) such that
\begin{align*}
&a_1 \not= a_2, \quad b_1 \not=b_2,\\
&\Phi^{-1}(\{a_i\}\times \C)\not= \emptyset\qquad (i=1,2),\\
&\Phi^{-1}(\C \times \{b_j\})\not= \emptyset\qquad (j=1,2).
\end{align*}
Let $\alpha_i \in \Phi^{-1}(\{a_i\}\times \C)$ ($i=1,2$) and
$\beta_j \in \Phi^{-1}(\C \times \{b_j\}$ ($j=1,2$).
By making use of the affine coordinate $(x,y)$ of $\C^2 \subset \pone^2$
we consider the following biholomorphic transform of $\pone^2$:
$$
\psi (x,y) = \left(\frac{x-a_1}{x-a_2}, \frac{y-b_1}{y-b_2}\right).
$$
Set $\Psi=\psi \circ \Phi: \C^2 \to \pone^2$.
Then points $\mathrm{P}_j$ ($1 \leq j \leq 4$) are exterior points of the
image $\Psi(\C^2)$.
Let $g:\C \to \C^2$ be a holomorphic curve such that the image
$g(\C)$ is contained by no analytic proper subset of $\C^2$,
and that $g$ passes through all four points $\alpha_i, \beta_j$
($i,j=1,2$).
Set $f=\Psi \circ g: \C \to \pone^2$.
Then $f$ is algebraically nondegenerate and $f(\C) \cap
E\not=\emptyset$.
By Proposition \ref{deg} $f$ is $\nabla$-nondegenerate, too.
\bigskip

{\it Problems.}  (i) It is an interesting problem to find more examples for
Theorem \ref{smt-1}.

(ii) It is naturally interesting to extend the results of \S5 to
the higher dimensional case and the case of general semi-abelian
varieties.

(iii) Is it possible to decrease ``2'' to ``1''
in the inequalities obtained by Theorems \ref{smt-5} and \ref{smt-7}.

\bigskip
\baselineskip=12pt
\rightline{Graduate School of Mathematical Sciences}
\rightline{The University of Tokyo}
\rightline{Komaba, Meguro,Tokyo 153-8914}
\rightline{e-mail: noguchi@ms.u-tokyo.ac.jp}

\begin{thebibliography}{99}
\setlength{\itemsep}{-3pt}
\bibitem{c33}
Cartan, H.,
Sur les z\'eros des combinaisons lin\'eaires de $p$ fonctions holomorphes donn\'ees,
Mathematica {\bf 7} (1933), 5-31.
\bibitem{f22}
Fatou, P.
Sur les fonctions m\'eromorphes de deux variables,
C. R. Acad.\ Sci.\ Paris {\bf 175} (1922), 862--865.
\bibitem{n77}
Noguchi, J.,
Holomorphic curves in algebraic varieties,
Hiroshima Math.\ J. {\bf 7} (1977), 833-853.
\bibitem{n81}
Noguchi, J.,
Lemma on logarithmic derivatives and holomorphic curves in algebraic varieties,
Nagoya Math.\ J.\ {\bf 83} (1981), 213-233.
\bibitem{n98}
Noguchi, J.,
On holomorphic curves in semi-Abelian varieties,
Math.\ Z. {\bf 228} (1998), 713-721.
\bibitem{n04}
Noguchi, J.,
Intersection multiplicities of holomorphic and algebraic curves with divisors,
Proc.\ OKA 100 Conference Kyoto/Nara 2001,
Advanced Studies in Pure Mathematics {\bf 42}, pp.\ 243-248,
Japan Math. Soc.\, Tokyo, 2004.
\bibitem{n08}
Noguchi, J.,
Value Distribution and Distribution of Rational Points
at Mittag-Leffler, 
Talk at Mittag-Leffler Institute, 27 March 2008.
\bibitem{no84}
Noguchi, J. and Ochiai, T.,
Geometric Function Theory in Several Complex Variables,
Japanese edition, Iwanami, Tokyo, 1984;
English Translation, Transl.\ Math.\ Mono.\ {\bf 80},
Amer.\ Math.\ Soc., Providence, Rhode Island,
1990.
\bibitem{nwy02}
Noguchi, J., Winkelmann, J. and Yamanoi, K.,
The second main theorem for holomorphic curves into semi-Abelian
varieties,
Acta Math.\ {\bf 188} no.1 (2002), 129-161.
\bibitem{nwy08}
Noguchi, J., Winkelmann, J. and Yamanoi, K.,
The second main theorem for holomorphic curves into semi-Abelian
varieties II,  Forum Math.{\bf 20} (2008), 469-503.
\bibitem{siu87}
Siu, Y.-T.,
Defect relations for holomorphic maps between spaces of different dimensions,
Duke Math.\ J. {\bf 55} (1987), 213-251.
\bibitem{t10}
Tiba, Y., Holomorphic curves into the product space of the Riemann
	spheres,
preprint UTMS 2010-19, 2010.
\bibitem{y04}
Yamanoi, K.,
Holomorphic curves in abelian varieties and intersection
with higher codimensional subvarieties, 
Forum Math.\ {\bf 16} (2004), {749}-{788}.
\end{thebibliography}
\end{document}